\RequirePackage{fix-cm}
\RequirePackage{amsmath}

\documentclass{article}
\usepackage{natbib}
\usepackage{arxiv}
\usepackage[english]{babel}
\usepackage{graphicx}
\usepackage{subfigure}
\usepackage{caption}

\usepackage[utf8]{inputenc} 
\usepackage[T1]{fontenc}    
\usepackage{hyperref}       
\usepackage{url}            
\usepackage{booktabs}       
\usepackage{amsmath,amsfonts}
\usepackage{amssymb}
\usepackage{nicefrac}       
\usepackage{microtype}      
\usepackage{lipsum}		
\usepackage{graphicx}
\usepackage{natbib}
\usepackage{doi}

\usepackage{hyperref}
\usepackage[table]{xcolor}
\usepackage{float}

\usepackage{textcomp}
\usepackage{gensymb}
\usepackage{booktabs}   
\usepackage{diagbox}

\DeclareMathOperator*{\argmin}{argmin}
\usepackage[T1]{fontenc}    
\usepackage[utf8]{inputenc} 
\usepackage[ruled,vlined]{algorithm2e}

\newtheorem{definition}{Definition}

\newcommand{\R}{\mathbb{R}}
\definecolor{LightCyan}{rgb}{0.88,1,1}
\definecolor{Gray}{gray}{0.9}

\title{Time Series Forecasting Using Manifold Learning}

\author{{Panagiotis Papaioannou}\\
	 Dipartimento di Matematica e Applicazioni\\
Universit\`a degli Studi di Napoli Federico II\\
	Naples, Italy\\
	\texttt{panagiotis.papaioannou@unina.it} \\
	\And
	{Ronen Talmon} \\
	Department of Electrical Engineering\\
	Technion — Israel Institute of Technology\\
	Haifa, Israel\\
	\texttt{ronen@ee.technion.ac.il} \\
	\AND
	{Ioannis G. Kevrekidis} \\
	Department of Chemical and Biomolecular Engineering \& \\Department of Applied Mathematics and Statistics \&\\
	Department of Medicine
	Johns Hopkins University\\
	Baltimore, USA\\
	\texttt{ronen@ee.technion.ac.il} \\
	\And
	{Constantinos Siettos} \thanks{Corresponding author} \\
 Dipartimento di Matematica e Applicazioni \&\\
 Scuola Superiore Meridionale\\
Universit\`a degli Studi di Napoli Federico II\\
	Naples, Italy\\
	\texttt{constantinos.siettos@unina.it} \\
	
}

\renewcommand{\shorttitle}{\textit{arXiv} Forecasting on Manifold}

\hypersetup{
pdftitle={Time Series Forecasting Using Manifold Learning},
pdfsubject={Numerical Analysis},
pdfauthor={Panagiotis Papaioannou, Ronen Talmon, Ioannis Kevrekidis Constantinos Siettos},
pdfkeywords={First keyword, Second keyword, More},
}

\begin{document}
\maketitle

\begin{abstract}
We address a three-tier numerical framework based on manifold learning for the forecasting of high-dimensional time series. At the first step, we embed the time series into a reduced low-dimensional space using a nonlinear manifold learning algorithm such as Locally Linear Embedding and Diffusion Maps. At the second step, we construct reduced-order regression models on the manifold, in particular Multivariate Autoregressive (MVAR) and Gaussian Process Regression (GPR) models, to forecast the embedded dynamics. At the final step, we lift the embedded time series back to the original high-dimensional space using Radial Basis Functions interpolation and Geometric Harmonics. For our illustrations, we test the forecasting performance of the proposed numerical scheme with four sets of time series: three synthetic stochastic ones resembling EEG signals produced from linear and nonlinear stochastic models with different model orders, and one real-world data set containing daily time series of 10 key foreign exchange rates (FOREX) spanning the time period 03/09/2001-29/10/2020. The forecasting performance of the proposed numerical scheme is assessed using the combinations of manifold learning, modelling and lifting approaches. We also provide a comparison with the Principal Component Analysis algorithm as well as with the naive random walk model and the MVAR and GPR models trained and implemented directly in the high-dimensional space.
\end{abstract}


\section{Introduction\label{Intro}}


Classical methods of forecasting include exponential smoothing and moving average models, ARIMA, Bayesian non parametric models including Gaussian Processes and other machine learning schemes such as Forward and Recurrent Neural Networks, Deep Learning, LSTMs, Reinforcement Learning, Reservoir Computing and Fuzzy Systems (\cite{Granger1969,de200625,lukovsevivcius2009reservoir,brockwell2016introduction,greff2016lstm,pathak2018model,vlachas2020backpropagation,wyffels2010comparative}). 

An inherent problem of using such data-driven modelling/forecasting approaches, when dealing with temporal measurements in a high-dimensional feature space, is the ``curse of dimensionality'' (\cite{bellmann1957dynamic}). In such cases, training, i.e.~the estimation of the values of the model parameters, requires a large number of snapshots, which increases exponentially with the dimension of the feature space. Thus, a fundamental task in modelling and forecasting is that of dimensionality reduction. Out of the many features that can be measured in time, one has first to identify the intrinsic dimension of the possible low-dimensional subspace (manifold) and the corresponding vectors that span it, which actually govern the emergent / macroscopically observed dynamics of the system under study. Assuming that the emergent dynamics evolve on a smooth enough low-dimensional manifold, an arsenal of linear, such as Singular value Decomposition (SVD) (\cite{golub1971singular}) and Dynamic Mode Decomposition (DMD) (\cite{schmid2010dynamic,mann2016dynamic,kutz2016dynamic}), and nonlinear such as kernel-PCA (\cite{scholkopf1998nonlinear}), Locally Linear Embedding (LLE) (\cite{roweis2000nonlinear}), ISOMAP (\cite{tenenbaum2000global,balasubramanian2002isomap}), Laplacian Eigenmaps (\cite{belkin2003laplacian}), Diffusion Maps (DMs) (\cite{coifman2005geometric,coifman2006diffusion,coifman2006geometric,nadler2006diffusion,coifman2008diffusion}), and Koopman operator (\cite{mezic2013analysis,williams2015data,dietrich2020koopman}) related manifold learning algorithms can be exploited towards this direction. 

For the two-fold task of dimensionality reduction and modelling in the low-dimensional subspace, \cite{chiavazzo2014reduced} constructed reduced kinetics models, by extracting the slow dynamics on a manifold globally parametrized by a truncated DM. A comparison of the reconstructed and the original high-dimensional dynamics was also reported. The reconstruction was achieved using Radial Basis Functions (RBFs) interpolation and Geometric Harmonics (GHs). \cite{liu2014coarse} used DMs to identify coarse variables that govern the emergent dynamics of a particle-based model of animal swarming and, based on these, they constructed a reduced stochastic differential equation. \cite{dsilva2015data} used DMs with a Mahalanobis distance as a metric to parametrize the slow manifold of multiscale dynamical systems. \cite{williams2015data} addressed an extension of DMD to compute approximations of the Koopman eigenvalues, thus showing that for large data sets the procedure converges to the numerical approximation one would obtain from a Galerkin method. The performance of the method was tested via the unforced Duffing equation and a stochastic differential equation with a double-well potential that 
converges to a Fokker-Planck equation. \cite{brunton2016discovering} used SVD to embed the dynamics of nonlinear systems in a linear manifold spanned by a few SVD modes, and constructed state-space linear models on the manifold to approximate the full dynamics. The so-called Sparse Identification of the Nonlinear Dynamics (SINDy) method was demonstrated using the Lorenz equations and the 2D Navier–Stokes equations.  \cite{bhattacharjee2016nonlinear} used ISOMAP to construct reduced-order models of heterogeneous hyperelastic materials in the 3D space, whose solutions are obtained by finite elements, while for the construction of the inverse map used RBF interpolation. \cite{wan2017reduced}  proposed a methodology for forecasting and quantifying
uncertainty in a reduced-dimensionality space, based on Gaussian process regression. The efficiency of the proposed approach was validated using data series from the Lorenz 96 model, the
Kuramoto-Sivashinsky, as well as a barotropic climate model in the form of a PDE. \cite{nathan2018applied} proposed a scheme based on the Koopman-operator theory to embed spatio-temporal dynamics of PDEs with the aid of DMD, for approximating the evolution of the high-dimensional dynamics on the reduced manifold; the reconstruction of the states of the original system was achieved by a linear transformation. \cite{chen2018molecular} have proposed a molecular enhanced sampling method based on the concept of autoencoders ~\cite{kramer1991nonlinear} to discover the low-dimensional projection of Molecular Dynamics and then reconstruct back the atomic coordinates.
\cite{vlachas2018data} proposed an LSTM-based method to predict the high-dimensional dynamics from the information acquired in a low-dimensional space. The embedding space is constructed based on the Takens embedding theorem (\cite{takens1981detecting}). The method is demonstrated through time-series  obtained from the Lorenz 96 equations, the Kuramoto–Sivashinsky equation and a barotropic climate model as in \cite{wan2017reduced}.~\cite{lee2020coarse} used DMs and Automatic Relevance Determination to construct embedded PDEs by the aid of Artificial Neural Networks (ANNs) and Gaussian Processes. The methodology was illustrated with Lattice-Boltzman simulations of the 1D Fitzhugh-Nagumo model. \cite{koronaki2020data} used DMs to embed the dynamics of a one-dimensional tubular reactor as modelled by a system of two PDEs on a 2D manifold, and then constructed a feedforward ANN to learn the dynamics on the 2D manifold. The efficiency of the scheme was compared with the original high-dimensional PDE dynamics by lifting the reduced dynamics back to the original space using GHs and RBFs. Recently, \cite{lin2021data} used the Koopman and the Mori-Zwanzig projection operator formalism to construct reduced-order models for chaotic and randomly forced dynamical systems, and then based on the Wiener projection they derived Nonlinear Auto-Regressive Moving Average with exogenous input models. The approach was demonstrated through the Kuramoto-Sivashinsky PDE and the viscous Burgers equation with stochastic forcing.

The performance of most of the above schemes was assessed in terms of interpolation, namely, the data are produced using dynamical models (ODEs, PDEs, SDEs or agent-based models) that were simulated in some specific interval of the high-dimensional domain, then reduced order models were constructed and trained based on the generated data, and finally a reconstruction of the high-dimensional dynamics was implemented within the original interval of the high-dimensional domain.

In this work, we assess the performance of a manifold learning numerical scheme within a different context: forecasting, i.e. out-of-sample extrapolation. This is a conceptually different task with respect to the interpolation problem of the model reduction of well-defined dynamical models in the form of ODEs or PDEs, as attempted in the above-mentioned works (see also the discussion in the conclusion section). In particular, We address a three-tier scheme to perform forecasting in the high-dimensional space, by (i) reducing the high-dimensional time series into a low-dimensional manifold using manifold learning, (ii) constructing and training reduced-order models on the manifold, and, based on these, making predictions, (iii) lifting the predictions in the high-dimensional space.
The performance of this ``embed-forecast-lift'' scheme is assessed through four sets of time series: three of them are synthetic time series resembling EEG recordings generated by linear and nonlinear discrete stochastic models with different model orders, while the fourth set is a real-world data set containing the time series of 10 key pairs of foreign exchange prices retrieved from the \url{www.investing.com} open API, spanning the period
03/09/2001-29/10/2020. For our computations, we used two well-established and widely applied nonlinear manifold learning algorithms, namely the LLE algorithm and DMs, two types of models, namely Multivariate Autoregressive (MVAR) and Gaussian Process Regression (GPR) models, and two operators for lifting, namely RBFs and GHs. We considered various combinations of the previous methods and compared them against the naive random walk and MVAR and GPR models implemented directly in the original space. A comparison with the results obtained with the leading PCA components for the case of the FOREX forecasting problem is also provided. To the best of our knowledge, this is the first work that addresses the problem of forecasting based on manifold learning, providing a comparison of various embedding, modelling and reconstruction approaches, and assessing their performance on both synthetic time series and a real-world FOREX data set.

\section{Methodology and preliminaries\label{MethodologyPrelimin}}

Let us denote by $\boldsymbol{x}_i \in \R^D$ the vector containing the features at the $i$-th snapshot of the time series, and by $\boldsymbol{X} \in \R^{D \times N}$ the matrix having as columns the vectors spanning a time window of $N$ observations. The assumption is that the data lie on a ``smooth enough'' low-dimensional (say of dimension $d$) manifold that is embedded in the high-dimensional $\R^D$ space. Our aim is to exploit manifold learning algorithms to forecast the time series in the high-dimensional feature space. Thus, our purpose is to bypass the ``curse of dimensionality'' associated with both the dimension of the original data and, importantly, with the limited size of their snapshots in time that usually characterize real-world data (such as financial data). Towards this aim, we propose a numerical scheme that consists of three steps.

At the first step, we employ embedding algorithms (LLE and DMs) to construct nonlinear maps from the high-dimensional to a low-dimensional subspace that preserve as much as possible the intrinsic geometry of the manifold (i.e.~maps assuring that neighborhoods of points in the high-dimensional space are mapped to neighborhoods of the same points in the low-dimensional manifold, and the notion of distance between the points is maintained as much as possible). In general it is expected that the manifold learning algorithms will provide an approximation to the intrinsic embedding dimension, which will differ from the ``true'' one. This depends on the threshold that one puts for the selection of the number of embedded vectors that span the low-dimensional manifold.
A central concept here is the Riemannian metric, which defines the properties of this map. At the second step, based on the resulting embedding features that span the low-dimensional subspace, we train a class of regression models (MVAR and GPR) to predict the evolution of the embedded time-series on the manifold based on their past values. The final step implements the lifting map. The aim here is to form the inverse map from the out-of-sample forecasted embedded points to the reconstruction of the features in the high-dimensional space. On one hand, the existence of such an inverse map is theoretically guaranteed by the assumption that a low-dimensional ``sufficiently smooth'' manifold exists. On the other hand, the data-driven derivation of the corresponding inverse map is neither unique (\cite{lee2013smooth}) nor trivial as for example when using PCA. In general, one has to solve a nonlinear least squares problem requiring the minimization of an objective function that can be formed using different criteria for the properties of the neighborhood of points; this results in different inverse maps whose performance has to be validated.

In what follows, for the completeness and clarity of the presentation, we very briefly present basic concepts of the manifold theory, that are useful to understand the steps of the proposed numerical methodology.

Manifold learning techniques can be viewed as unsupervised learning algorithms in the sense that they ``learn'' from the available data a low-dimensional representation of the original high-dimensional space, thus providing an ``optimal'' (under certain assumptions) embedded subspace where the information available in the high-dimensional space is preserved
as much as possible. Here, we briefly present some basic elements of the theory of manifolds and manifold learning (see e.g. \cite{berger2012differential,kuhnel2015differential,lee2013smooth,wang2012geometric}).
Let us start with the definition of a $d$-dimensional manifold.
\begin{definition}[Manifold]
\textit{A set $M \subset \R^{n}$ is called a $d$-dimensional manifold (of class $C^{\infty}$) if for each point $\boldsymbol{p} \in M$, there is an open set $W \subset M$ such that there exists a $C^{\infty}$ bijection $\boldsymbol{f}: U \rightarrow W$, $U \subset \R^d$ open, with $C^{\infty}$ inverse $\boldsymbol{f}^{-1}: W \rightarrow U$ (i.e. $W$ is $C^{\infty}$-diffeomorphic to $U$). The bijection $\boldsymbol{f}$ is called a local parametrization of $M$, while $\boldsymbol{f}^{-1}$ is called a local coordinate mapping, and the pair $(W,\boldsymbol{f}^{-1})$ is called a chart (or neighborhood) of $\boldsymbol{p}$ on $M$.}
\end{definition}

\noindent
Thus, in a coordinate system $(W,\boldsymbol{f}^{-1})$, a point $\boldsymbol{p} \in W$ can be expressed by the coordinates $(f^{-1}_1,f^{-1}_2,\dots, f^{-1}_n)$, where $f^{-1}_i$ is the $i$-th element of $\boldsymbol{f}^{-1}$. A chart $(W,\boldsymbol{f}^{-1})$ is centered at $\boldsymbol{p}$ if and only if $\boldsymbol{f}^{-1}(\boldsymbol{p}) = \boldsymbol{0}$.
Thus, we always assume that the manifold satisfies the $T_2$ Hausdorff separation axiom, stating that every pair of distinct points on $M$ have disjoint open neighborhoods. Summarizing, a manifold is a Hausdorff (Separated/T2) space, in which every point has a neighborhood that is diffeomorphic to an open subset in $\R^d$.

\begin{definition}[Tangent space and tangent bundle to a manifold]
\textit{Let $\boldsymbol{p}\in M \subset \R^n$. The tangent space to $M$ at $\boldsymbol{p}$ is the vector subspace $T_{\boldsymbol{p}}M$ formed by the tangent vectors at $\boldsymbol{p}$ defined by
$$ T_{\boldsymbol{p}}M = D\boldsymbol{f}_{\boldsymbol{u}}(T_{\boldsymbol{u}}\R^{d}), \quad \boldsymbol{f}(\boldsymbol{u})=\boldsymbol{p}, $$
where $T_{\boldsymbol{u}}\R^{d}$
is the tangent space at $\boldsymbol{u}$ and
\begin{equation*}
    D\boldsymbol{f}_{\boldsymbol{u}}(T_{\boldsymbol{u}}\R^{d})=\left[\frac{\partial f_i(\boldsymbol{u})}{\partial u_j} \right]=
    \begin{bmatrix} \frac{\partial f_1}{\partial u_1} & \frac{\partial f_1}{\partial u_2} & \dots & \frac{\partial f_1}{\partial u_d}\\
    \frac{\partial f_2}{\partial u_1} & \frac{\partial f_2}{\partial u_2} & \dots & \frac{\partial f_2}{\partial u_d}\\
    \vdots & \vdots & \vdots & \vdots\\
    \frac{\partial f_n}{\partial u_1} & \frac{\partial f_n}{\partial u_2} & \dots & \frac{\partial f_n}{\partial u_d}
    \end{bmatrix},
\end{equation*}
with $rank(D\boldsymbol{f})=d$. Thus, $S= \left\{ \frac{\partial \boldsymbol{f}}{\partial u_1}, \frac{\partial \boldsymbol{f}}{\partial u_2}, \dots, \frac{\partial\boldsymbol{f}}{\partial u_d} \right\}$ is a basis for  $T_{\boldsymbol{p}}M$.
The union of tangent spaces
$$TM = \cup_{\boldsymbol{p} \in M} T_{\boldsymbol{p}}M$$
is called the tangent bundle of $M$.} 
\end{definition}


\begin{definition}[Riemannian manifold and metric] 
\textit{A Riemannian manifold is a manifold endowed with a positive definite inner product $g_{\boldsymbol{p}}$ defined on the tangent space $T_{\boldsymbol{p}}M$ at each point $\boldsymbol{p}$. The family of inner products $g_{\boldsymbol{p}}$ is called a Riemannian metric and the Riemannian manifold is denoted $(M,g)$.}
\end{definition}


The above definitions refer to the case of continuum limit -- infinite number of points. However, for real data with a limited number of observations, one can only try to approximate the manifold; the analytic charts of the Riemannian metric that actually define the geometry of the manifold are simply not available. Thus, a consistent manifold learning algorithm constructs in a numerical way a map that approximates in a probabilistic way the continuous one when the sample size $N \to \infty$. Thus, a fundamental preprocessing step of a class of manifold learning algorithms such as Laplacian Maps, LLE, ISOMAP and DMs is the construction of a weighted graph, say $\mathcal{G}(\boldsymbol{X},E)$, where $E$ denotes the set of weights between points. This construction is based on a predefined metric (such as the Gaussian kernel and the $k$-NN algorithm), which is used to appropriately connect each point $\boldsymbol{x}\in\R^D$ with all the others. This defines a weighted graph of neighborhoods of points. Theoretically, with an appropriate choice of the metric and its parameters, the graph is guaranteed to be connected, i.e.~there is a path between every pair of points in $\boldsymbol{X}$ (\cite{lee2013smooth}).

For the completeness of the presentation, in the following sections we briefly discuss the manifold learning algorithms, the regression models and the lifting techniques used in this work.

\section{Manifold learning algorithms\label{ManifoldLearning}}

\subsection{Locally Linear Embedding (LLE)\label{LocalLinearEmbeddingFramework}}

The LLE algorithm (\cite{Saul2003}) is a nonlinear manifold learning algorithm which constructs $\mathcal{G}(\boldsymbol{X},E)$ based on the $k$-NN algorithm. It is assumed that the neighborhood forms a basis for the reconstruction of any point in the neighborhood itself. Thus, every point is written as a linear combination of its neighbors. The weights of all pairs are then estimated with least squares, minimizing the $L_2$ norm of the difference between the points and their neighborhood-based reconstruction. The same rationale is assumed for the low-dimensional subspace, keeping the same estimates of the weights in the high-dimensional space. In the low-dimensional subspace, one now seeks for the coordinates of the points that minimize the $L_2$ norm of the difference between the coordinate of the points on a $d-$dimensional manifold and the weighted neighborhood reconstruction. The minimization problem is represented by an eigenvector problem, whose first $d$ bottom non-zero eigenvectors provide an ordered set of orthogonal coordinates that span the manifold. It should be noted that in order to obtain a unique solution to the minimization problem, the number of $k$ nearest neighbors should not exceed the dimension of the input space (\cite{Saul2003}).

Thus, the LLE algorithm can be summarized in the following three basic steps (\cite{Saul2003}).
\begin{enumerate}
\item Identify the $k$ nearest neighbors for all $\boldsymbol{{x}_{i}} \in R^D , \; i=1,2,\dots,N$ (with $k\ge d + 1$).
   For each $\boldsymbol{{x}}_{i}$ this forms a set $K\{\boldsymbol{{x}}_{i}\} \subset \R^{k \times D}$ containing the $k$ nearest neighbors of $\boldsymbol{{x}}_{i}$. 
\item  Write each $\boldsymbol{{x}}_{i}$ in terms of $K\{\boldsymbol{{x}}_{i}\}$ as 
   \begin{equation}
   \boldsymbol{{x}}_{i} = \sum_{j \in K\{\boldsymbol{{x}}_{i}\}} w_{ij}\boldsymbol{{x}}_{j},
   \label{eq1LLE}
   \end{equation}
   and find the matrix $\boldsymbol{W} = [w_{ij}] \in \R^{N \times N}$ of the unknown weights by minimizing the objective function
   \begin{equation}
   \mathcal{L}(\boldsymbol{W}) = \sum_{i=1}^{N} \left\lVert \boldsymbol{{x}}_{i} - \sum_{j=1, j\neq i}^{N} w_{ij}\boldsymbol{x}_{j} \right\rVert^{2}_{L_2}, \label{eq2LLE}
   \end{equation}
   with the constraint
   \begin{equation}
   \sum_{j=1}^{N} w_{ij} = 1.
   \end{equation}
\item Embed the points $\boldsymbol{{x}}_{i} \in \R^D$, $i=1,2,\dots, N$, into a low-dimensional space with coordinates $\boldsymbol{y}_{i} \in \R^d$, $i=1,2,\dots, N$, $d \ll D$. This step in LLE is accomplished by computing  the vectors $\boldsymbol{y}_{i} \in \R^d$ that minimize the objective function:
   \begin{equation}
   \phi(\boldsymbol{Y})=\sum_{i=1}^{N} \left\lVert \boldsymbol{y}_i - \sum_{j=1, j\neq i}^{N} w_{ij}\boldsymbol{y}_{j} \right\rVert ^{2}_{L_2}
   \label{eq3LLE},
   \end{equation}
   where the weights $w_{ij}$ are fixed at the values found by solving the minimization problem \eqref{eq2LLE}. The embedding vectors are required to be centered at the origin with an identity covariance matrix (\cite{Saul2003}).  The embedded vector $\boldsymbol{y}_i$ are constrained so that they have a zero mean and a unit covariance matrix.
\end{enumerate}

The cost function \eqref{eq3LLE} is quadratic and can be stated as
\begin{equation}
\phi(\boldsymbol{Y})=\sum_{i=1}^{N} \sum_{j=1}^{N} Q_{ij}\langle \boldsymbol{y}_{i}, \boldsymbol{y}_{j} \rangle,
\end{equation}
involving the inner products of the embedding vectors and a square matrix $\boldsymbol{Q}$ with elements
\begin{equation}
Q_{ij} = \delta_{ij} - w_{ij} - w_{ji} + \sum_k w_{ki} w_{kj},
\end{equation}
$\delta_{ij}=1$ if $i=j$, and $0$ otherwise. In matrix form, $\boldsymbol{Q}$ can be written as an $N \times N$ \emph{sparse} symmetric matrix:
\begin{equation}
\boldsymbol{Q} = (\boldsymbol{I}-\boldsymbol{W})^T (\boldsymbol{I}-\boldsymbol{W}).
\end{equation}
In practice, this sparse representation of $\boldsymbol{Q}$ gives rise to a significant computational reduction, especially when $N$ is large, since left multiplication by $\boldsymbol{Q}$ is given by
\begin{equation}
\boldsymbol{Q}\boldsymbol{v} = (\boldsymbol{v}-\boldsymbol{W}\boldsymbol{v}) - \boldsymbol{W}^T(\boldsymbol{v}-\boldsymbol{W}\boldsymbol{v}), \quad \boldsymbol{v} \in \mathbb{R}^N,
\end{equation}
requiring just two multiplications by the sparse matrices $\boldsymbol{W}$ and $\boldsymbol{W}^T$.

Minimizing the cost function \eqref{eq3LLE} can be computed via the eigenvalue decomposition of $\boldsymbol{Q}$. Specifically, the eigenvector associated with the smallest (zero) eigenvalue is the vector with  all 1's  and it trivially minimizes \eqref{eq3LLE} (\cite{Saul2003}); it is disregarded because it leads to a constant (degenerate) embedding. The optimal embedding is therefore obtained by the $d$ eigenvectors of $\boldsymbol{M}$, denoted as $\boldsymbol{q}_k \in R^N$, $k=1,2,\dots,d$, corresponding to the next $d$ smallest eigenvalues. Thus, the coordinates of $\boldsymbol{x}_i$, $i=1,2,\dots,N$, in the embedded space are given by the vector
\begin{equation}
\boldsymbol{\mathcal{R}}(\boldsymbol{x}_i)=\boldsymbol{y}_i=[q_{i1},q_{i2},\dots,q_{id}]^T,
\end{equation}
where $q_{ij}$ denotes the $j$-th element of eigenvector $\boldsymbol{q}_{i}$.

\subsection{Diffusion Maps (DMs)\label{DiffusionMaps}}
Here, the construction of the affinity matrix is based on the computation of a random walk on the graph $\mathcal{G}(\boldsymbol{X},E)$. The first step is to construct a graph using a kernel function,
say $k(\boldsymbol{x}_i,\boldsymbol{x}_j)$.
The kernel function can be chosen as a Riemannian metric, so that it is symmetric and positive definite. 
Standard kernels, such as the Gaussian kernel, typically define a neighborhood of each point $\boldsymbol{x}_i$, i.e.~a set of points $\boldsymbol{x}_j$ having strong connections with $\boldsymbol{x}_i$, namely, large values of $k(\boldsymbol{x}_i,\boldsymbol{x}_j)$.

At the next step, one constructs a Markovian transition matrix, say $\boldsymbol{P}$, whose elements correspond to the probability of jumping from one point to another in the high-dimensional space.
%
This transition matrix defines a Markovian (i.e.~memoryless) random walk $X_t$ by setting
$$ p_{ij} = p(\boldsymbol{x}_i,\boldsymbol{x}_j) = \text{Prob}(X_{t+1} = \boldsymbol{x}_j|X_t = \boldsymbol{x}_i).$$ 

For a graph constructed from a sample of finite size $N$, the weighted degree of a point (node) is defined by
\begin{equation}
deg(\boldsymbol{x}_i) = \sum_{j=1}^N k(\boldsymbol{x}_i,\boldsymbol{x}_j),
\end{equation}
and the volume of the graph is given by
$$ vol(\mathcal{G}) = \sum_{i=1}^N d(\boldsymbol{x}_i).$$
Then, the random walk on such a weighted graph can be defined by the transition probabilities
\begin{equation}
p_{ij} = p(\boldsymbol{x}_i,\boldsymbol{x}_j) = \frac{k(\boldsymbol{x}_i,\boldsymbol{x}_j)}{deg(\boldsymbol{x}_i)}.
\label{transprob}
\end{equation}
Clearly, from the above definition, we have that $\sum_{j=1}^{N} p(\boldsymbol{x}_i,\boldsymbol{x}_j)=1$. 

We note that in the continuum, the above can be described as a continuous Markov process on a probability space $(\Omega, \mathcal{H}, \mathcal{P})$, where $\Omega$ is the sample space, $\mathcal{H}$ is a $\sigma$-algebra of events in $\Omega$, and $\mathcal{P}$ is a probability measure. 
Let $\mu$ be the density function of the points in the sample space $\Omega$ induced from the probability measure, $\mu: \mathcal{H} \rightarrow \R$ with $\mu(\Omega)=1$.
Then, using the kernel function $k$, the transition probability from a point $\boldsymbol{x}\in \Omega$ to another point $\boldsymbol{y} \in \Omega$ is given by
\begin{equation}
p(\boldsymbol{x},\boldsymbol{y})=\frac{k(\boldsymbol{x},\boldsymbol{y})}{
\int_{\boldsymbol{\Omega}} k(\boldsymbol{x},\boldsymbol{y}) d\mu(\boldsymbol{\Omega})},
\end{equation}
which is the continuous-space counterpart of \eqref{transprob}.

The above procedure defines a row-stochastic transition matrix, $\boldsymbol{P} = [p_{ij}]$, which encapsulates the information about the neighborhoods of the points. 
Note that by raising $\boldsymbol{P}$ to the power of $t=1,2,\dots $, we get the jumping
probabilities in $t$ steps. This way, the underlying geometry is revealed through high or low transition probabilities between the points, i.e. paths that follow the underlying geometry have a high probability of occurrence, while paths away from the ``true'' embedded structure have a low probability. Note that the $t$-step transition probabilities, say $p_t(\boldsymbol{x}_i,\boldsymbol{x}_j)$, satisfy the Chapman-Kolmogorov equation: 
\begin{equation}
    p_{t_1+t_2}(\boldsymbol{x}_i,\boldsymbol{x}_j)=\sum_{\boldsymbol{x}_k\in \boldsymbol{X}}p_{t_1}(\boldsymbol{x}_i,\boldsymbol{x}_k)p_{t_2}(\boldsymbol{x}_k,\boldsymbol{x}_j).
\end{equation}
The Markov process defined by the probability matrix $\boldsymbol{P}$ has a stationary distribution given by
\begin{equation}
    \pi(\boldsymbol{x}_i)=\frac{deg(\boldsymbol{x}_i)}{\sum_{\boldsymbol{x}_j\in \boldsymbol{X}}deg(\boldsymbol{x}_j)},
\end{equation}
and it is reversible, i.e.
\begin{equation}
    \pi(\boldsymbol{x}_i)p(\boldsymbol{x}_i,\boldsymbol{x}_j)= \pi(\boldsymbol{x}_j)p(\boldsymbol{x}_j,\boldsymbol{x}_i), \quad \forall \boldsymbol{x}_i, \boldsymbol{x}_i \in \boldsymbol{X}.
\end{equation}

Furthermore, if the kernel is appropriately chosen so that the graph is connected, then the Markov chain is ergodic and irreducible. The Brouwer Fixed Point Theorem (\cite{kellogg1976constructive}) implies that the transition matrix of an ergodic process has a stationary vector $\boldsymbol{\pi}$ such that:
\begin{equation}
    \boldsymbol{P}^T\boldsymbol{\pi}=\boldsymbol{\pi},
\end{equation}
and hence the matrix of a Markov ergodic chain has always an eigenvalue 1 corresponding to the stationary state. From the Perron-Frobenius Theorem, we know that its geometric multiplicity is 1. It can be shown that all other eigenvalues have a magnitude smaller than 1 (\cite{coifman2006diffusion}).
From \eqref{transprob} we get
\begin{equation}
\boldsymbol{P}=\boldsymbol{D}^{-1}\boldsymbol{K}, \quad \boldsymbol{D}=\text{diag} \left( \sum_{j=1}^{N}k_{ij} \right),
\label{diagonP}
\end{equation}
where $k_{ij} = k(\boldsymbol{x}_i,\boldsymbol{x}_j)$.
We note that the transition matrix $\boldsymbol{P}$ is similar to the symmetric and positive definite matrix $\boldsymbol{S}=\boldsymbol{D}^{-1/2}\boldsymbol{K}\boldsymbol{D}^{-1/2}$. Thus, the transition matrix $\boldsymbol{P}$ has a decomposition given by
\begin{equation}
\boldsymbol{P}=\sum_{i=1}^N \lambda_i \boldsymbol{u}_i \boldsymbol{v}_{i}^{T}, 
\end{equation}
where $\lambda_i \in \R$ are the (positive) eigenvalues of $\boldsymbol{P}$, $\boldsymbol{u}_i \in \R^N$ are the left eigenvectors, and $\boldsymbol{v}_i \in \R^N$ are the right eigenvectors, such that $\langle \boldsymbol{u}_i,\boldsymbol{v}_{j} \rangle =
\delta_{ij}$.

The set of right eigenvectors $\boldsymbol{v}_i$ establishes an orthonormal basis for the subspace $R(\boldsymbol{P}^T)$ of $\R^D$ spanned by the rows of $\boldsymbol{P}$. Row $i$ represents the transition probabilities from point $\boldsymbol{x}_i$ to all the other points of the graph. According to the Eckart-Young-Mirsky Theorem (\cite{eckart1936approximation,mirsky1960symmetric}), the best $d$-dimensional low-rank approximation of the row space of $\boldsymbol{P}$ in the Euclidean space $\R^d$ is provided by its $d$ right eigenvectors corresponding to the $d$ largest eigenvalues.

It has been shown that asymptotically, when the number of data points uniformly sampled from a low dimensional manifold goes to infinity, the matrix $\frac{1}{\sigma}(\boldsymbol{I}-\boldsymbol{P})$ (where $\sigma$ is the kernel function scale) approaches the Laplace-Beltrami operator of the underlying Riemannian manifold (\cite{coifman2005geometric,nadler2006diffusion}). This allows to consider the eigenvectors corresponding to the largest few eigenvalues of $\boldsymbol{P}$ as discrete approximations of the principal eigenfunctions of the Laplace-Beltrami operator. Since the principal eigenfunctions of the Laplace-Beltrami operator establish an accurate embedding of the manifold (\cite{jones2008manifold}), this result promotes the usage of the eigenvectors of $\boldsymbol{P}$ for practical (discrete) data embedding.

At this point, the so-called diffusion distance, which is an affinity distance related to the reachability between two points $\boldsymbol{x}_i$ and $\boldsymbol{x}_j$, is given by
\begin{equation}
D_{t}^2(\boldsymbol{x}_i,\boldsymbol{x}_j) = ||p_t(\boldsymbol{x}_i,\cdot),p_t(\boldsymbol{x}_j,\cdot)||_{L_2,1/deg}^2 = \sum_{i=1}^{N} \frac{(p_t(\boldsymbol{x}_i,\boldsymbol{x}_k)-p_t(\boldsymbol{x}_j,\boldsymbol{x}_k))^2} {deg(\boldsymbol{x}_k)},
\label{diffdist}
\end{equation}
where $p_t(\boldsymbol{x}_i,\cdot)$ is the $i$-th row of $\boldsymbol{P}^t$.
The embedding of the $t$-step transition probabilities is achieved by forming a family of maps (DMs) of the $N$ points $\boldsymbol{x}_i \in \boldsymbol{X}$ in the Euclidean subspace of $\R^d$ defined by
\begin{equation}
 \boldsymbol{\mathcal{R}}_t(\boldsymbol{x}_i)=\boldsymbol{y}_i=\left[\lambda_{1}^t\boldsymbol{v}_1(i),\lambda_{2}^t \boldsymbol{v}_2(i),\dots, \lambda_{d}^t \boldsymbol{v}_d(i)\right]^T, i=1,2,\dots N.
\end{equation}
A useful property of DMs is that the Euclidean distance in the embedded space $||\boldsymbol{\mathcal{R}}_t(\boldsymbol{x}_i)-\boldsymbol{\mathcal{R}}_t(\boldsymbol{x}_j)||_{L_2}$ is the best $d$-dimensional approximation of the diffusion distance, given by \eqref{diffdist}. 
where the equality holds for $d=N$.


\section{Forecasting with Regression Models\label{ForecastingModels}}

Let us assume that the time series are being generated by a (nonlinear) law of the form
\begin{equation}
    \boldsymbol{y}_{t}=\boldsymbol{\phi}(\boldsymbol{z},\boldsymbol{\mu})+\boldsymbol{e}_t,
    \label{regressmodel}
\end{equation}
where $\boldsymbol{y}_{t} \in \R^d$ denotes the vector containing the measured values of the response variables at time $t$, $\boldsymbol{\phi}: \R^p \times \R^q \rightarrow \R^d $ is a smooth function encapsulating the law that governs the system dynamics, which depends on parameters represented by the vector $\mu \in \R^q$, $\boldsymbol{z} \in \R^p$ is the vector containing the explanatory input variables (predictors), which may include past values of the response variables at certain times, say $t-1, t-2, \dots, t-m$, and other exogenous variables, and $\boldsymbol{e}_t$ is the vector of the unobserved noise at time $t$.

In general, the forecasting problem (here in the low-dimensional space) using regression models can be written as a minimization problem of the following form:
\begin{equation}
    \argmin_{\boldsymbol{g} \in G, \; \boldsymbol{\theta} \in \R^{l}}{\mathcal{L}(\boldsymbol{y}_{t+k}-\boldsymbol{g}(\boldsymbol{\boldsymbol{z},\boldsymbol{\theta}}))},
\end{equation}
where $k$ is the prediction horizon, $\boldsymbol{g}: \R^p \times \R^l \rightarrow \R^d$ is the regression model
\begin{equation}
\boldsymbol{\hat{y}}_{t+k}=\boldsymbol{g}(\boldsymbol{\boldsymbol{z},\boldsymbol{\theta}}),
\end{equation}
with parameters
$\boldsymbol{\theta} \in \R^l$, $G$ is the space of available models $\boldsymbol{g}$, $\hat{\boldsymbol{y}}_{t}$ is the output of the model at time $t$, and $\mathcal{L}$ is the loss function that determines the optimal forecast.

We note that the forecasting problem can be posed in two different modes (\cite{marcellino2006comparison}): the iterative and the direct one. In the iterative mode, one trains one-step ahead models based on the available time series and simulates/iterates the models to predict future values. In the direct mode, forecasting is performed in a sliding-window framework, where the model is retrained with the data contained within the sliding window to provide a multiperiod-ahead value of the dependent variables. Here, we aim at testing the performance of the proposed scheme for one-period ahead of time predictions (i.e. with $k=1$) with both iterative and directs modes. For our illustrations, we used two forecasting models, namely MVAR and GPR. A brief description of the models follows.

\subsection{Multivariate Autoregressive (MVAR) model\label{MVARsec}}

An MVAR model can then be written as
\begin{equation}
\boldsymbol{y}_t=\boldsymbol{\theta_0}+\sum_{j=1}^{m}\boldsymbol{y}_{t-j}\boldsymbol{\Theta_{j}}+\boldsymbol{e}_t,
\label{MVAR}
\end{equation}
where $\boldsymbol{y}_t = [y_{t1}, y_{t2}, \dots, y_{td}]^T \in \R^d$ is the vector of the response time series at time $t$, $m$ is the model order, i.e. the maximum time lag,  $\boldsymbol{\theta_0} \in \R^d$ is the regression intercept, $\boldsymbol{\Theta}_{j} \in \R^{d \times d}$ is the matrix containing the regression coefficients of the MVAR model, and $\boldsymbol{e}_t=[e_{t}^{(1)}, e_{t}^{(2)}, \dots, e_{t}^{(d)}]^T$ is the vector of the unobserved errors at time $t$, which are assumed to be uncorrelated random variables with zero mean and constant variance $\sigma^2$.
In a more general form, in view of \eqref{regressmodel}, the MVAR model can be written as
\begin{equation}
y_{ik}=\theta_{0k}+\sum_{j=1}^{m}\theta_{jk}z_{ij}+e_{ik}, \quad i=1,2,\dots d, \quad k=1,2,\dots N,
\label{mvarcovfullindex}
\end{equation}
where $y_{ik}$ is the model output for the $i$-th variable at the $k$-th time instant, $\theta_{0k} \in \R$ is the corresponding regression intercept and $\theta_{jk}$ the corresponding $j$-th regression coefficient, $z_{ij}$ is the $j$-th predictor of the $i$-th response variable (e.g. the time-delayed time series). 
According to the Gauss-Markov theorem, the best unbiased linear estimator of the regression coefficients is the one that results from the solution of the least-squares (LS) problem
\begin{equation*}
\arg\min_{\theta_{0k}, \theta_{jk}} \sum_{i=1}^{d} \sum_{k=1}^{N}
\left( y_{ik} - \theta_{0k} - \sum_{j=1}^{m} \theta_{jk}z_{ij} - e_{ik} \right)^2 ,
\end{equation*}
given by
\begin{equation}
    \boldsymbol{\hat{\Theta}}=(\boldsymbol{Z}^T\boldsymbol{Z})^{-1}\boldsymbol{Z}^T\boldsymbol{Y},
\end{equation}
where $\boldsymbol{Y} = [\boldsymbol{y}_1,\boldsymbol{y}_2,\dots \boldsymbol{y}_d] \in R^{N \times d}$ and $\boldsymbol{Z}=[\boldsymbol{1}_N, \boldsymbol{z}_1,\boldsymbol{z}_2,\dots,\boldsymbol{z}_m] \in \R^{N\times (m+1)}$.
Assuming that the unobserved errors are i.i.d. normal random variables, then by the maximum likelihood estimate of the error covariance, one can also estimate the forecasting intervals of a new observation.

\subsection{Gaussian Process Regression (GPR)\label{GPR}}
For the implementation of GPR it is assumed that the unknown function $\boldsymbol{\phi}$ in \eqref{regressmodel} can be modelled by $d$ single-output Gaussian distributions given by
\begin{equation}
    P(\boldsymbol{\phi}_i|\boldsymbol{z})=\mathcal{N}(\boldsymbol{\phi}_i|\boldsymbol{\mu},\boldsymbol{K}(\boldsymbol{z}, \boldsymbol{z} |\boldsymbol{\theta})),
\end{equation}
where $\boldsymbol{\phi}_i=[\phi_i(\boldsymbol{z}_1), \phi_i(\boldsymbol{z}_2), \dots, \phi_i(\boldsymbol{z}_N)]$, and $\phi_i$ is the $i$-th component of $\boldsymbol{\phi}$, $\boldsymbol{\mu}(\boldsymbol{z})$ is the vector with the expected values of the function, and $\boldsymbol{K}(\boldsymbol{z}, \boldsymbol{\theta})$ is a $N\times N$ covariance matrix formed by a kernel. The prior mean function is often set to $\boldsymbol{\mu}(\boldsymbol{z}) = \boldsymbol{0}$ with appropriate normalization of the data.

Predictions at a new point, say $\boldsymbol{z}_{*}$, are made by drawing $\boldsymbol{\phi}_{i*}$ from the posterior distribution $P(\boldsymbol{\phi}_i |(\boldsymbol{Z},\boldsymbol{y}_i))$, where $\boldsymbol{Z}=[\boldsymbol{z}_1,\boldsymbol{z}_2,\dots,\boldsymbol{z}_N]^T \in \R^{N \times p}$, and by assuming a Gaussian distributed noise, and appropriate normalization of the data points is given by the joint multivariate normal distribution:
\begin{equation}
    \begin{bmatrix}
    \boldsymbol{y}_i\\
    \phi_i(\boldsymbol{x}_*)
    \end{bmatrix}\sim \mathcal{N}
    \begin{pmatrix} 
    \boldsymbol{0}\\
    0
    \end{pmatrix},
    \quad
    \begin{bmatrix}
    \boldsymbol{K}(\boldsymbol{Z},\boldsymbol{Z}|\boldsymbol{\theta})+\sigma^2\boldsymbol{I} & \boldsymbol{k}(\boldsymbol{Z},\boldsymbol{z}_*\boldsymbol{\theta}) \\
     \boldsymbol{k}(\boldsymbol{z}_*,\boldsymbol{Z}|\boldsymbol{\theta}) & k(\boldsymbol{z}_*,\boldsymbol{z}_*|\boldsymbol{\theta})
    \end{bmatrix}, \quad i=1,2,\dots,d.
\end{equation}
It can be shown that the posterior conditional distribution
\begin{equation}
    P(\phi_i(\boldsymbol{x}_*)|\boldsymbol{y}_i,\boldsymbol{Z},\boldsymbol{z}_*)
\end{equation}
can be analytically derived and the expected value and covariance of the estimation are given by:
\begin{equation}
    {\boldsymbol{\bar{\phi}}_i}(\boldsymbol{z}_*)=\boldsymbol{k}(\boldsymbol{z}_*,\boldsymbol{Z}|\boldsymbol{\theta})
    \left( \boldsymbol{K}(\boldsymbol{Z},\boldsymbol{Z}|\boldsymbol{\theta})+\sigma^2\boldsymbol{I} \right)^{-1} \boldsymbol{y}_i, \quad i=1,2,\dots, d
    \label{meangp}
\end{equation}
\begin{equation}
    \sigma_{*}^2={k}(\boldsymbol{z}_*,\boldsymbol{z}_*|\boldsymbol{\theta})-\boldsymbol{k}(\boldsymbol{z}_*,\boldsymbol{Z}|\boldsymbol{\theta})[ \boldsymbol{K}(\boldsymbol{Z},\boldsymbol{Z}|\boldsymbol{\theta})+\sigma^2\boldsymbol{I}]^{-1}\boldsymbol{k}(\boldsymbol{Z},\boldsymbol{z}_*|\boldsymbol{\theta})
    \label{covgp}
\end{equation}
The hyperparameters in the above equations are estimated by minimizing
the marginal likelihood that is given by
\begin{align*}
    \ell & = -\log p(\boldsymbol{y}_i|\boldsymbol{Z},\boldsymbol{\theta}) \\
         & = \frac{1}{2}\boldsymbol{y}_{i}^T[ \boldsymbol{K}(\boldsymbol{Z},\boldsymbol{Z}|\boldsymbol{\theta})+\sigma^2\boldsymbol{I}]^{-1}\boldsymbol{y}_{i}+\frac{1}{2}\log[\boldsymbol{K}(\boldsymbol{Z},\boldsymbol{Z}|\boldsymbol{\theta})+\sigma^2\boldsymbol{I}]+\frac{N}{2}\log 2\pi .
\end{align*}
Here, for our computations we have used a mixed kernel composed by adding the fundamental components for forecasting with GPR (see e.g. \cite{corani2021time}), namely a radial basis function kernel:
\begin{equation}
    k(\boldsymbol{z}_i,\boldsymbol{z}_j)=\theta_{1}^2 \exp \left( -\frac{||\boldsymbol{z}_i-\boldsymbol{z}_j||^{2}_{L_2}}{2\theta_{2}^2} \right),
\end{equation}
a linear kernel:
\begin{equation}
    k(\boldsymbol{z}_i,\boldsymbol{z}_j)=\theta_{3}^2+\theta_{4}^2\langle \boldsymbol{z}_i,\boldsymbol{z}_j\rangle,
\end{equation}
a periodic kernel:
\begin{equation}
    k(\boldsymbol{z}_i,\boldsymbol{z}_j) = \theta_{5}^2 \exp \left( -\frac{2\sin^2{(\pi||\boldsymbol{z}_i-\boldsymbol{z}_j||_{L_2}/\theta_6)}}{\theta_{7}^2} \right),
\end{equation}
and a white noise kernel:
\begin{equation}
    k(\boldsymbol{z}_i,\boldsymbol{z}_j)=\theta_{8}^2\delta_{ij} .
\end{equation}

\section{Reconstruction of the high-dimensional space: solving the pre-image problem\label{LiftingOperators}}

The final step is to reconstruct the high-dimensional space from measurements on the manifold, i.e., to ``lift''  the predictions made by the reduced-order models on the manifold back to the original high-dimensional space. In the case of PCA, this task is trivial, as the solution to the reconstruction problem, given by
\begin{equation}
    \argmin_{\boldsymbol{U}_d\in \R^{D\times d}}\sum_{i=1}^N||\boldsymbol{x}_i-\boldsymbol{U}_d\boldsymbol{U}_{d}^{T}\boldsymbol{x}_i||_{L_2}^2,  \quad \boldsymbol{U}_{d}^{T}\boldsymbol{U}_d=I,
\end{equation}
is just a linear transformation which maximizes the variance of the data on the linear manifold, and is given by the first $d$ principal eigenvectors of the covariance matrix.

In the case of nonlinear manifold learning algorithms, such as LLE and DMs, we want to learn the inverse map (lifting operator):
\begin{equation}
    \boldsymbol{\mathcal{L}} 
    \equiv \boldsymbol{\mathcal{R}}^{-1}: \boldsymbol{\mathcal{R}}(\boldsymbol{X}) \rightarrow \boldsymbol{X},
\end{equation}
for new samples on the manifold $\boldsymbol{y}_* \not\in \boldsymbol{\mathcal{R}}(\boldsymbol{X})$. This inverse problem is referred as ``out-of-sample extension pre-image'' problem. The ``out-of-sample extension'' problem usually refers to the direct problem, i.e. that of learning the direct embedding map (i.e. the restrictions to the manifold) $\boldsymbol{\mathcal{R}}(\boldsymbol{X}): \boldsymbol{X} \rightarrow \boldsymbol{\mathcal{R}}$, for new samples in the input space $\boldsymbol{x}_* \not\in 
\boldsymbol{X}$. Towards this aim, a well established  methodology is the  Nystr\"{o}m extension (\cite{nystrom1929praktische}).

In general, the solution of the pre-image problem can be posed as:
\begin{equation}
    \arg\min_{\boldsymbol{c}} ||\boldsymbol{y}_*-\boldsymbol{\mathcal{R}}(\boldsymbol{\mathcal{L}}(\boldsymbol{y}_*)|\boldsymbol{c}))|| ,
    \label{preimage}
\end{equation}
subject to a constraint, where $\boldsymbol{\mathcal{L(\cdot)|\boldsymbol{c}}}$ is the lifting operator depending on some parameters $\boldsymbol{c}$.\par

Below, we describe the reconstruction methods that we used in this work, namely, RBFs interpolation and Geometric Harmonics (GH), which provide a direct solution of the inverse problem, thus giving some insight about their implementation and pros and cons. For a review and a comparison of such methods in the framework of chemical kinetics see~\cite{chiavazzo2014reduced}.

\subsection{Radial Basis Function (RBF) Interpolation\label{RBF}}

The lifting operator is constructed with interpolation through RBFs among the corresponding set of neighbors of the new state $\boldsymbol{y}_*$ in the ambient space. The lifting operator is defined by (\cite{chiavazzo2014reduced}):
\begin{equation}
    \boldsymbol{\mathcal{L}}(\boldsymbol{y}_*) = x_{i*} = \sum_{j=1}^k c_{ji} \psi(||\boldsymbol{y}_*-\boldsymbol{y}_j||), \quad i=1,2,\dots, D ,
    \label{liftrbf}
\end{equation}
where $x_{i*}$ is the $i$-th coordinate of $\boldsymbol{x}_*$, the $\boldsymbol{y}_j$'s are the neighbors of the unseen sample $\boldsymbol{y}_*$ on the manifold, and $\psi$ is the radial basis function.

Similarly, the restriction operator can be written as:
\begin{equation}
    \boldsymbol{\mathcal{R}}(\boldsymbol{x}_*) =y_{i*} = \sum_{j=1}^k c_{ji}\psi(||\boldsymbol{\mathcal{L}}(\boldsymbol{y}_*)-\boldsymbol{\mathcal{L}}(\boldsymbol{y}_j)||), \quad i=1,2,\dots,D,
    \label{restrictrbf}
\end{equation}
where $\boldsymbol{\mathcal{L}}(\boldsymbol{y}_j)=\boldsymbol{x}_j$ is the known image of $\boldsymbol{y}_j$ (the neighbors of $\boldsymbol{y}_*$ in the ambient space).
The unknown coefficients of the lifting operator can be computed by setting at the left hand side of \eqref{liftrbf} the coordinates of the one-to-one corresponding neighbors of $\boldsymbol{y}_*$ in the ambient space. Here, we consider a special class of RBFs known as radial powers, given by $$\psi(||\boldsymbol{\mathcal{L}}(\boldsymbol{y}_*)-\boldsymbol{\mathcal{L}}(\boldsymbol{y}_j)||)=||\boldsymbol{\mathcal{L}}(\boldsymbol{y}_*)-\boldsymbol{\mathcal{L}}(\boldsymbol{y}_j)||^{p}_{L_2},$$
where $p$ is an odd integer. By doing so, the unknown coefficients $c_{ji}$ of the lifting operator are given by the solution of the following linear system:
\begin{equation}
\boldsymbol{A}
   \begin{bmatrix}
   c_{1i}\\
   c_{2i}\\
   \vdots\\
   c_{ki}\end{bmatrix}=
    \begin{bmatrix}
   x_{1i}\\
   x_{2i}\\
   \vdots\\
   x_{ki}
   \end{bmatrix}, \quad i=1,2,\dots,D,
   \label{linsys_rbf}
\end{equation}
where 
\begin{equation}
\boldsymbol{A }= \begin{bmatrix}
   0& ||\boldsymbol{y}_1-\boldsymbol{y}_2||^{p}_{L_2}&\dots & ||\boldsymbol{y}_1-\boldsymbol{y}_{k}||^{p}_{L_2}\\
   ||\boldsymbol{y}_2-\boldsymbol{y}_1||^{p}_{L_2}& 0&\dots & ||\boldsymbol{y}_2-\boldsymbol{y}_{k}||^{p}_{L_2}\\
   \vdots & \vdots & \dots & \vdots\\
   ||\boldsymbol{y}_{k}-\boldsymbol{y}_1||^{p}_{L_2}& ||\boldsymbol{y}_{k}-\boldsymbol{y}_2||^{p}_{L_2}&\dots & 0 \end{bmatrix},
   \label{solverbf}
\end{equation}
where $x_{ji}$ is the $i$-th coordinate of the $j$-th point in the ambient space, whose restriction on the manifold is the $j$-th nearest neighbor of $\boldsymbol{y}_*$. For our computations, we used $p=1$.
Then, \eqref{liftrbf} can be used to find the coordinates of $\boldsymbol{x}_*$ in the ambient space.

We should note that radial powers are better for the task of interpolation compared to Gaussian kernels, as they do not suffer from ill conditioning. The matrix $\boldsymbol{A}$ is invertible with the assumption that the centers $\boldsymbol{y}_j$ are distinct (see e.g. \cite{monnig2014inverting,amorim2015facing}). However, $\boldsymbol{A}$ can be rank deficient from a numerical point of view, e.g. because the points may be very close to each other. Thus, the solution of \eqref{solverbf} using a method such as the LU factorization of $\boldsymbol{A}$ may result in numerical instabilities.

\subsection{Geometric Harmonics (GHs)\label{GH}}

GHs are a set of functions that allow the extension of the embedding of new unseen points on the manifold $\boldsymbol{x}_* \not\in \boldsymbol{X}$, which are not given in the set of points used for building the embedding (\cite{coifman2006geometric}). Their derivation is based on the Nystr\"{o}m (or quadratic) extension method (\cite{nystrom1929praktische}), which has been used for the 
numerical solutions of integral equations (\cite{delves1988computational,press1990fredholm}), and in particular the Fredholm equation of second kind, reading:
\begin{equation}
    f(t)=g(t)+\mu \int_{a}^{b} k(t,s)f(s) ds,
    \label{fredholm}
\end{equation}
where $f(t)$ is the unknown function, while $k(t,s)$ and $g(t)$ are known. The Nystr\"{o}m method starts with the approximation of the integral, i.e.
\begin{equation}
    \int_{a}^{b}y(s)ds \approx \sum_{j=1}^N w_j y(s_j),
    \label{quadint}
\end{equation}
where $s_j$ are $N$ appropriately chosen collocation points and $w_j$ are the corresponding weights, which are determined, e.g. by the Gauss-Jacobi quadrature rule. Then, by using \eqref{quadint} in \eqref{fredholm} and evaluating $f$ and $g$ at the $N$ collocation points, we get the following approximation:
\begin{equation}
    (\boldsymbol{I}-\mu \boldsymbol{\tilde{K}})\boldsymbol{\hat{f}}=\boldsymbol{g},
\end{equation}
where the matrix $\boldsymbol{\tilde{K}}$ has elements  $\tilde{k}_{ij}  = k(s_i,s_j) w_j$.
Based on the above, the solution of the homogenous Fredholm problem ($\boldsymbol{g}=0$) is given by the solution of the eigenvalue problem
\begin{equation}
   \boldsymbol{\tilde{K}}\boldsymbol{\hat{f}}=\frac{1}{\mu}\boldsymbol{\hat{f}},
\end{equation}
i.e.
\begin{equation}
    \sum_{j=1}^N w_j k(s_i,s_j) \hat{f}_j = \frac{1}{\mu}\hat{f}_i, \quad i=1,2,\dots, N,
\end{equation}
where $\hat{f}_i = \hat{f}(s_i)$ is the $i$-th component of $\boldsymbol{\hat{f}}$. The Nystr\"{o}m extension of $f(t)$, using a set of $N$ sample (collocation) points,
at an arbitrary point $x$ in the full domain is given by:
\begin{equation}
   \boldsymbol{\mathcal{E}}({f}(x)) = \hat{f}(x) = \mu \sum_{j=1}^N w_j k(x,s_j) \hat{f}_j.
\end{equation}
Within the framework of DMs, we seek for the out-of-the-sample (filtered) extension of a real-valued function $f$ defined at $N$ sample points $\boldsymbol{x}_i \in \boldsymbol{X}$ to one or more unseen points $\boldsymbol{x}_* \notin \boldsymbol{X}$. The function $f$ can be for example a DM coordinate $\lambda_{j}^t \boldsymbol{v}_j(\boldsymbol{x}_i)$, $j=1,2,\dots d$, or another function representing the output of a regression model (see also the discussion in \cite{thiem2020emergent,rabin2016earthquake}.

Recalling that the eigenvectors $\boldsymbol{v}_j$ form a basis, the extension is implemented (see \cite{coifman2006geometric}) by first expanding $f(\boldsymbol{x}_i)$ in the first $d$ parsimonious eigenvectors $\boldsymbol{v}_l$ of the Markovian matrix $\boldsymbol{P}^t$:
\begin{equation*}
 \hat{f}(\boldsymbol{x}_i) = \sum_{l=1}^d a_l \boldsymbol{v}_l(\boldsymbol{x}_i), \quad i=1,2,\dots N,
\end{equation*}
where $a_l=\langle \boldsymbol{v}_l, \boldsymbol{f}\rangle$ are the projection coefficients of the function on the first $d$ parsimonious eigenvectors, and $\boldsymbol{f} \in \R^N$ is the vector containing the values of the function $f$ at the $N$ points $\boldsymbol{x}_i$.
Then, one computes the Nystr\"{o}m extension of $\hat{f}$ at $\boldsymbol{x}_*$ using the same projection coefficients as
\begin{equation}
 \boldsymbol{\mathcal{E}}(\hat{f}(\boldsymbol{x}_*)) = \sum_{l=1}^d a_l \boldsymbol{\hat{v}}_l(\boldsymbol{x}_*),
\end{equation}
where 
\begin{equation}
   \boldsymbol{\hat{v}}_l(\boldsymbol{x}_*)=\frac{1}{\lambda_{l}^t}\sum_{j=1}^N k(\boldsymbol{x}_j,\boldsymbol{x}_*) \boldsymbol{v}_{l}(\boldsymbol{x}_j), \quad l=1,2,\dots, d.
\end{equation}
are the corresponding GHs. 
Scaling up the (filtered) extension of the function $f$ to a set of say $L$ new points can be computed using the following matrix product (\cite{coifman2006geometric,thiem2020emergent}):
\begin{equation}
   \boldsymbol{\mathcal{E}}(\hat{\boldsymbol{f}})=\boldsymbol{K}_{L\times N}\boldsymbol{V}_{N\times d}{\boldsymbol{\Lambda}_{d\times d}^{-1}}\boldsymbol{V}_{d\times N}^T\boldsymbol{f}_{N\times 1},
\end{equation}
where $\boldsymbol{K}_{L\times N}$ is the corresponding kernel matrix, $\boldsymbol{V}_{N\times d}$ is the matrix with columns the $d$ parsimonious eigenvectors $\boldsymbol{v}_l$, $\boldsymbol{\Lambda}_{d\times d}$ is the diagonal matrix with elements $\lambda_{l}^t$.

The above direct approach provides a map from the ambient space to the reduced-order space (restriction) and vice versa (lifting).

\section{The forecasting problems\label{Forecasting}}

For demonstrating the performance of the proposed numerical framework and comparing the various embedding, modelling and reconstruction approaches, we used three synthetic time series resembling EEG recordings produced by linear and nonlinear stochastic discrete models, and a real-world FOREX pair data set spanning the period from 03/09/2001 until 29/10/2020. Below, we describe the models and the FOREX data set.

\subsection{The synthetic time series\label{Thesynthetictimeseries}}

Our synthetic stochastic signals resemble EEG time series (see e.g. \cite{baccala2001partial,nicolaou2016nonlinear}) produced by linear and nonlinear five-dimensional discrete stochastic models with white noise.

The linear five-dimensional stochastic discrete model is given by the following equations:
\begin{eqnarray}\label{linearmodel}
\nonumber
y^{(1)}_t=0.2y^{(1)}_{t-1}-0.4y^{(2)}_{t-1}+w^{(1)}_t ,
\\ \nonumber
y^{(2)}_t=-0.5y^{(1)}_{t-1}+0.15y^{(2)}_{t-1}+w^{(2)}_t ,
\\ 
y^{(3)}_t=-0.14y^{(2)}_{t-1}+w^{(3)}_t ,
\\ \nonumber 
y^{(4)}_t=0.5y^{(1)}_{t-1}-0.25y^{(2)}_{t-1}+w^{(4)}_t ,
\\ \nonumber
y^{(5)}_t=0.15y^{(1)}_{t-1}+w^{(5)}_t,
\end{eqnarray}
where for the training process the model order is assumed to be known (here equal to 1).

A second (nonlinear) model is given by the following equations (see e.g. \cite{nicolaou2016nonlinear}):
\begin{eqnarray}\label{nonlinearmodel}
\nonumber
y^{(1)}_t=3.4y^{(1)}_{t-1}(1-{y^{(1)}_{t-1}}^2)\exp{(-{y^{(1)}_{t-1}}^2)}+w^{(1)}_t ,
\\ \nonumber
y^{(2)}_t=3.4y^{(2)}_{t-1}(1-{y^{(2)}_{t-1}}^2)\exp{(-{y^{(2)}_{t-1}}^2)}+0.5y^{(1)}_{t-1}y^{(2)}_{t-1}+w^{(2)}_t ,
\\ 
y^{(3)}_t=3.4y^{(3)}_{t-1}(1-{y^{(3)}_{t-1}}^2)\exp{(-{y^{(3)}_{t-1}}^2)}+0.3y^{(2)}_{t-1}+0.5{y^{(1)}_{t-1}}^2+w^{(3)}_t ,
\\ \nonumber 
y^{(4)}_t=0.5y^{(1)}_{t-1}-0.25y^{(2)}_{t-1}+w^{(4)}_t ,
\\ \nonumber
y^{(5)}_t=0.15y^{(1)}_{t-1}+w^{(5)}_t.
\end{eqnarray}

Finally, the proposed scheme was also validated through the time series produced by a linear stochastic model with a model order greater than 1 given by the following equations: 
\begin{eqnarray}\label{linearmodel3}
\nonumber
y^{(1)}_t=0.1y^{(1)}_{t-1}-0.6y^{(2)}_{t-3}+w^{(1)}_t ,
\\ \nonumber
y^{(2)}_t=-0.15y^{(1)}_{t-3}+0.8y^{(2)}_{t-3}+w^{(2)}_t ,
\\ 
y^{(3)}_t=-0.45y^{(2)}_{t-3}+w^{(3)}_t ,
\\ \nonumber 
y^{(4)}_t=0.45y^{(1)}_{t-3}-0.85y^{(2)}_{t-3}+w^{(4)}_t ,
\\ \nonumber
y^{(5)}_t=0.95y^{(1)}_{t-2}+w^{(5)}_t.
\end{eqnarray}

\noindent
In all the above models, the time series $w^{(i)}_t$, $i=1,2,3,4,5$, are uncorrelated normally-distributed noise with zero mean and unit standard deviation. 

\subsection{The FOREX forecasting problem\label{Headings}}

We used the \texttt{investpy 0.9.14} module of Python (\cite{Bartolome2020}) to download 10 USD-based FOREX pairs daily spot closing prices from the \url{www.investing.com} open API, spanning the period from 03/09/2001 until 29/10/2020:
EUR/USD, GBP/USD, AUD/USD, NZD/USD, JPY/USD, CAD/USD, CHF/USD, SEK/USD, NOK/USD, DKK/USD.
The time period was selected to contain crucial market crashes like the 2008 Lehman related stock market crash, the 2010 sovereign debt crisis in the Eurozone and even the most recent Covid-19 related market crash of 2020.

For our analysis, we computed the compounded daily returns of the FOREX pairs as
\begin{equation}
{r}_{i,t} = \log \left( \frac{S_{i,t}}{S_{i,t-1}} \right),
\label{eq:FXRawReturns}
\end{equation}
where $S_{i,t}$ is the spot price of the $i$-th FOREX pair in the data set above, and $r_{i,t}$ is the corresponding compounded return at date $t$.
However, trading the FOREX markets is not only associated with the spot price movements but also with the so-called ``currency carry trading''. The currency carry-trading involves first the identification of high interest-paying investment assets (e.g. bonds or short-term bank deposits) denominated in a country currency (e.g. Japanese Yen). Then, the investment is carried on by borrowing money in another currency, for which the paying interest rate is lower. Such a trading requires the approximation of the so called interest-rate-differential (IRD) excess returns, which should be incorporated in the FOREX price movements (see e.g. \cite{menkhoff2012carry}). Here, we used the short term interest rates data retrieved from the OECD database as the proxy of these IRD excess returns. Since the short term interest rates data are reported on a monthly basis and on an annual percentage format, we constructed daily approximations, by linearly interpolating through the available downloaded records (using the \texttt{Pandas} module of Python (\cite{McKinney2010})) and normalizing on the basis of an annual calendar period. 

After this pre-processing step, we denote as IR$_{\text{USA},t}$ the time series of the daily approximations of the USA short term interest rates,
and IR$_{i,t}$ is the corresponding time series of each of the other $i=1,\dots, 10$ countries. Then, the interest rates differentials (IRD) time series are given by:
\begin{equation}
    \text{IRD}_{i,t} = \text{IR}_{i,t} - \text{IR}_{\text{USA},t}.
\label{eq:IRD}
\end{equation}
Finally, the so called ``carry adjusted returns'' of the FOREX are given by (see also e.g. \cite{menkhoff2012carry}):
\begin{equation}
x_{i,t} = r_{i,t} + \text{IRD}_{i,t},
\label{eq:FXCarryAdjustedReturns}
\end{equation}
where $x_{i,t}$ is the $i$-th FX pair ``carry'' adjusted return 
corresponding to its raw ``unadjusted'' market price return, $r_{i,t}$ is defined in \eqref{eq:FXRawReturns}, and IRD$_{i,t}$.

For quantifying the potential excess returns (profits), we constructed a trading strategy based on the so-called risk parity rationale (see e.g. \cite{Braga2015}) where each asset is allocated with a portfolio weight which is proportional to its inverse risk. Thus, one assigns higher portfolio weights to less volatile assets, and smaller portfolio weights to more risky assets. Thus, the risk is quantified using the volatility $\sigma_{i, t}$ (measured by the standard deviation of logarithmic returns over a specific time period up to the $t$ trading day) of each asset $i$ of the portfolio plus the total portfolio volatility. In our FOREX problem, the risk parity portfolio allocation practically means investing $1/\sigma_{i, t}$ at each of the $i=1,2,\dots,10$ FOREX pairs with corresponding carry-adjusted returns $x_{i,t}$. By performing one-step ahead predictions for each of carry-adjusted returns of the 10 pairs, denoted as $\hat{x}_{i,t+1}$, we create a ``binary'', or as otherwise called ``directional'', trading signal of ``buy'' or ``sell'' for each one of the FOREX. Thus, the  trading strategy reads as follows:
\begin{equation}
    u_{i,t} = sign(\hat{x}_{i,t+1}) = \\
    \begin{cases}
    \phantom{-}1, \; \text{``buy''} \\
    -1, \; \text{``sell''} \\
    \end{cases}
\end{equation}

Based on the above, the profit or loss at the next day ($t+1$) is given by:
\begin{equation}
\Pi_{t+1} = \sum_{i=1}^D \frac{u_{i,t} x_{i,t+1}}{{\sigma}_{i,t}},
\label{eq:RiskParityPortfolio}
\end{equation}
where $x_{i,t+1}$ is the real $i$-th FOREX pair return at time $t+1$, and $\Pi_{t+1}$ is the risk parity portfolio return at time $t+1$.

\section{Numerical Results}
\label{NumericalResults}

For the implementation of the numerical algorithms, we used the \texttt{datafold}, \texttt{sklearn} and \texttt{statsmodels} packages of Python (\cite{Lehmberg2020,scikit-learn,seabold2010statsmodels}). 
The selection of the eigensolver was based on the sparsity and size of the input matrix: ARPACK was used if the size of the input matrix was greater than 200 and $n+1<10$ (where $n$ is the number of requested eigenvalues), otherwise, the ``dense'' option was used. The ARPACK package (\cite{lehoucq1998arpack}) implements implicitly restarted Arnoldi methods, using a random default starting vector for each iteration with a tolerance of $tol=10^{-6}$ and a maximum number of 100 iterations. This approach is implemented by the function \texttt{scipy.sparse.linalg.eigsh} of the \texttt{scipy} module (upon which the \texttt{sklearn} one depends) (\cite{Virtanen2020}). The ``dense'' eigensolver is implemented by the function \texttt{eigh} of the \texttt{scipy} module and returns the eigenvalues and eigenvectors computed using the LAPACK routine \texttt{syevd} using the divide and conquer algorithm (\cite{Cuppen1980,lapack99}). We used the default value of the tolerance of the Newton-Raphson iterations, which is of the order of the floating-point precision, and the maximum number of iterations was set to $30N$ iterations, where $N$ is the size of the matrix.

The DM embedding on the training set of all data sets was performed using $d$ parsimonious eigenvectors (\cite{dsilva2018parsimonious}). Here, for the construction of the graph at the first step, we used the standard Gaussian kernel, defined by
\begin{equation}
    k(\boldsymbol{x}_i, \boldsymbol{x}_j) = \exp{\left( -\frac{\| \boldsymbol{x}_i - \boldsymbol{x}_j \|_{L_2}^2}{\sigma} \right)},
\end{equation}
where the $\| \boldsymbol{x}_i - \boldsymbol{x}_j \|_{L_2}$ is the Euclidean distance between $\boldsymbol{x}_i$ and $\boldsymbol{x}_j$, and $\sigma$ is a scaling parameter that controls the size of the neighborhood (or the connectivity of the graph). For its derivation, we follow the systematic approach provided by \cite{singer2009detecting}. The full kernel is used for the calculations. 

The VAR models were trained using the VAR class of the \texttt{statsmodels.tsa.} \texttt{vector\_ar.var\_model} routine, using the OLS default method for the parameters estimation. The corresponding LAPACK function used to solve the least-squares problem is the default \texttt{gelsd}, which exploits the QR or LQ factorization of the input matrix. 

The hyperparameters of the GPR model were optimized using the (default) L-BFGS-B algorithm of the \texttt{scipy.optimize.minimize} method (\cite{Byrd1995,Zhu1997}). The gradient vector was estimated using forward finite differences with the numerical approximations of the Jacobian being performed with a default step size $eps=10^{-8}$. The iterations were stopped either when 
$$ \frac{(f^k - f^{k+1})}{\max\{|f^k|,|f^{k+1}|,1\}} < 10^{-9}, $$
$f^k$ being the value of the loss function at step $k$, or when
$$ \max_i | \text{proj}(g)_i | < 10^{-5}, $$
$\text{proj}(g)_i$ being the $i$-th component of the projected gradient at the current iterate.
We used the default values for the maximum number of function evaluations and the maximum number of iterations (15,000) as well as the default value for the maximum number of line searches per iteration (20). 

For the lifting task, 50 nearest neighbors were considered for interpolation by all methods. The underlying $k$-NN algorithm is based on the algorithm proposed in \cite{Maneewongvatana2002}. Using different values of the number of nearest neighbors within the range 20-100 did not change the outcomes of the analysis. In the case of RBF interpolation, the underlying linear system of equations was solved by the LAPACK \texttt{dgesv} routine from \texttt{scipy.linalg.lapack.dgesv}, which implements the default method of the LU decomposition with partial pivoting. For the GH approach, we used the Gaussian kernel. In effect, we are performing ``double" DM here - computing diffusion maps on the leading retained diffusion map components for the reduced embedding. This procedure, suggested in ~\cite{chiavazzo2014reduced} can actually be performed ``once and for all" globally~\cite{Evangelou}.

The computations were performed using a system with an Intel Core i7-8750H CPU @2.20 GHz and 16GB of RAM. 

\subsection{Synthetic time series}

For both the linear and nonlinear models, we produced 2000 points. We used 1500 points for learning the manifold and for training the various models, and 500 points to test the forecasting performance of the various schemes. The forecasting performance was tested using the iterative mode, i.e. by training the models for one-step ahead predictions and then simulating the trained model iteratively to predict future values. The performance was measured using the Root Mean Square Error (RMSE) of the residuals, reading:
\begin{equation}
    \mbox{RMSE} = \frac{1}{N} \sqrt{\sum_{i=1}^N (\hat{\boldsymbol{x}_i} - \boldsymbol{x}_i) ^ 2},
    \label{MSE}
\end{equation}
where $\hat{\boldsymbol{x}_i}$ are the predictions and $\boldsymbol{x}_i$ the actual data. To quantify the forecasting intervals due to the stochasticity of the models, we performed 100 runs for each combination of manifold learning algorithms (DMs and LLE), models (MVAR and GPR) and lifting methods (RBFs and GHs), reporting the median and the 5-th and 95-th percentiles of the resulting RMSE. The RMSE values obtained with the naive random walk model, as well as with the MVAR and GPR models trained in the original space, are also given for comparison purposes.

In Table \ref{LinearManifoldTable}, we report the forecasting statistics of the time series produced with the linear model given by \eqref{linearmodel} as obtained over 100 runs. As it is shown, both the MVAR and GPR models trained in the original five-dimensional space outperform the naive random walk. The RMSEs of the MVAR and GPR models suggest a good match with the stochastic model, with the residuals being approximately within one standard deviation of the distribution of the noise level. In the same table, we provide the corresponding forecasting statistics as obtained with the proposed ``embed-forecast-lift'' scheme for the various combinations of manifold learning algorithms, regression models and lifting approaches.
\begin{table}[ht]
 \caption{Linear Stochastic Model~\eqref{linearmodel}. RMSE statistics (median, 5th and 95th percentiles over 100 runs) for each of the five variables as obtained by: (a) training MVAR and GPR models with model order one in the original 5D space (MVAR(OS), GPR(OS)), and (b) the proposed ``embed-forecast-lift'' scheme for all combinations of the manifold learning algorithms (DMs and LLE) with two coordinates, models (MVAR and GPR) and lifting approaches (RBFs and GHs). For comparison purposes the RMSEs obtained with the naive random walk is also reported.}
 \centering
 \addtolength{\tabcolsep}{-5pt}
  \begin{tabular}{|l|c|c|c|c|c|}
  \hline
  \rowcolor{Gray}
    Model/Variable & $y^{(1)}_t$ & $y^{(2)}_t$ & $y^{(3)}_t$ & $y^{(4)}_t$ & $y^{(5)}_t$\\
    \hline
    \hline
    \rowcolor{LightCyan}
    Random Walk & 1.374  & 1.446  & 1.427  & 1.551  & 1.425  \\
    \rowcolor{LightCyan}
  &  (1.287,1.466) &  (1.333,1.528) & (1.339,1.497) &  (1.455,1.654) & (1.335,1.512)\\
    \hline
    MVAR(OS) & 1.151  & 1.180  & 1.010  & 1.224  & 1.011  \\
    & (1.077,1.218) &  (1.113,1.269) &  (0.966,1.062) & (1.156,1.281) & (0.960,1.063) \\
    \hline
     \rowcolor{LightCyan}
    GPR (OS) & 1.151 & 1.179 & 1.010 & 1.224 & 1.011 \\
     \rowcolor{LightCyan}
    & (1.077,1.218) &  (1.113,1.272) & (0.966,1.061) &  (1.154,1.284) & (0.959,1.064) \\
    \hline
    \hline 
    DM-GPR-GH  & 1.153  & 1.184  & 1.012  & 1.228  & 1.013  \\
     &  (1.082,1.220) &  (1.117,1.273) & (0.966,1.061) & (1.160,1.287) & (0.959,1.065) \\
     \hline
      \rowcolor{LightCyan}
    LLE-GPR-GH & 1.155 & 1.182 & 1.011 & 1.230 & 1.014 \\
     \rowcolor{LightCyan}
    & (1.077,1.218) &  (1.115,1.272) &  (0.966,1.065) & (1.157,1.292) & (0.962,1.064) \\
    \hline
    \hline 
    DM-GPR-RBF & 1.260  & 1.390  & 1.163  & 1.365  & 1.166  \\
   & (1.111,2.201) & (1.149,1.980) &  (0.997,2.331) &  (1.173,1.939) & (0.997,1.853)\\
   \hline
    \rowcolor{LightCyan}
    LLE-GPR-RBF & 1.197 & 1.234 & 1.076 & 1.270 & 1.103  \\
     \rowcolor{LightCyan}
      & (1.103,1.452) & (1.131,1.461) & (0.993,1.582) & (1.175,1.468) & (0.990,1.553) \\
    \hline
    \hline 
    DM-MVAR-GH & 1.151  & 1.180  & 1.011  & 1.225  & 1.011  \\
      & (1.078,1.220) & (1.113,1.269) & (0.966,1.063) & (1.155,1.288) & (0.959,1.063) \\
      \hline
       \rowcolor{LightCyan}
    LLE-MVAR-GH & 1.155  & 1.182  & 1.011 & 1.229  & 1.014  \\
     \rowcolor{LightCyan}
     & (1.077,1.218) & (1.115,1.271) &  (0.966,1.065)  &  (1.158,1.293) &  (0.962,1.064) \\
    \hline
    \hline 
    DM-MVAR-RBF & 1.289  & 1.312  & 1.135  & 1.341  & 1.106  \\
    & (1.124,2.132) & (1.146,2.018) & (0.996,1.635) & (1.186,1.933) & (0.99,1.601) \\
    \hline
     \rowcolor{LightCyan}
    LLE-MVAR-RBF  & 1.187  & 1.230 & 1.080  & 1.264  & 1.096  \\
     \rowcolor{LightCyan}
   & (1.106,1.477) & (1.14,1.477) & (0.995,1.602) & (1.177,1.479) & (0.989,1.558) \\
    \hline
  \end{tabular}
  \label{LinearManifoldTable}
\end{table}
For our illustrations, we have chosen the first two parsimonious DM coordinates, and the corresponding two LLE eigenvectors. As shown, the best performance is obtained with the GH lifting operator. Using GHs for lifting and any combination of the selected manifold learning algorithms and models outperforms all other combinations, thus resulting in practically the same RMSE values when compared with the predictions made in the original space. This suggests that the proposed ``embed-forecast-lift'' scheme applied in the 5D feature space provides a very good reconstruction of the predictions made in the original space. On the other hand, lifting with RBF interpolation with LU decomposition, generally resulted in poor reconstructions of the high-dimensional space, thus, in many cases, giving wide forecasting intervals that contained the median value of the naive random walk RMSE.

Next, in Table \ref{NonLinearManifoldTable}, we report the forecasting statistics of the time series for the nonlinear stochastic model \eqref{nonlinearmodel} as obtained over 100 runs. As in the case of the linear stochastic model, both MVAR and GPR trained in the original 5D space outperform the naive random walk. The resulting RMSE values suggest a good match with the nonlinear stochastic model. Yet, the match is poorer than the one obtained for the linear model.
\begin{table}[ht]
 \caption{Nonlinear Stochastic Model \eqref{nonlinearmodel}. RMSE statistics (median, 5th and 95th percentiles over 100 runs) for each of the five variables as obtained by: (a) training MVAR and GPR in the original 5D space (MVAR(OS), GPR(OS)) with model order one, (b) using the proposed ``embed-forecast-lift'' scheme for all the combinations of the manifold learning algorithms (DMs and LLE), models (MVAR and GPR), and lifting approaches (RBFs and GHs). For the embedding, three coordinates were used.}
\centering
 \addtolength{\tabcolsep}{-5pt}
  \begin{tabular}{|l|c|c|c|c|c|}
  \hline
  \rowcolor{Gray}
 Model/Variable & $y^{(1)}_t$ & $y^{(2)}_t$ & $y^{(3)}_t$ & $y^{(4)}_t$ & $y^{(5)}_t$\\
 \hline
 \hline
  \rowcolor{LightCyan}
 Random Walk & 1.711  & 2.085  & 2.292  & 1.731  & 1.435 \\
  \rowcolor{LightCyan}
 & (1.597,1.815) & (1.939,2.258) & (2.147,2.437) & (1.62,1.819) & (1.354,1.527) \\
 \hline
 MVAR(OS) & 1.181  & 1.438  & 1.565  & 1.212  & 1.021   \\
  &  (1.123,1.234) & (1.359,1.555) & (1.463,1.648)& (1.151,1.272) & (0.967,1.076) \\
    \hline
     \rowcolor{LightCyan}
    GPR(OS) & 1.182 & 1.437 & 1.684  & 1.213 & 1.020 \\
     \rowcolor{LightCyan}
     & (1.123,1.233) & (1.360,1.555) & (1.540,1.799) & (1.151,1.272) & (0.966,1.076) \\
    \hline
 \hline 
    DM-GPR-GH  & 1.181  & 1.440  & 1.573  & 1.214  & 1.022 \\
     & (1.123,1.234) & (1.363,1.56) & (1.477,1.654) & (1.151,1.276) & (0.968,1.077)\\
     \hline 
      \rowcolor{LightCyan}
    LLE-GPR-GH & 1.182 & 1.451 & 1.579 & 1.214 & 1.023  \\
     \rowcolor{LightCyan}
      & (1.126,1.236) & (1.363,1.563) & (1.473,1.675) & (1.151,1.273) & (0.97,1.076) \\
    \hline
    \hline 
    DM-GPR-RBF  & 1.488  & 1.618  & 2.060  & 1.710 & 1.199 \\
      & (1.172,8.88) & (1.392,9.95) & (1.642,6.102) & (1.203,16.868) & (1.012,9.108)\\
    \hline
     \rowcolor{LightCyan}
    LLE-GPR-RBF & 1.214 & 1.449 & 1.587 & 1.254 & 1.088   \\
     \rowcolor{LightCyan}
       & (1.133,1.557) & (1.365,1.561) & (1.502,1.692) & (1.159,1.446) & (0.989,1.5) \\
    \hline
    \hline 
    DM-MVAR-GH & 1.181  & 1.438  & 1.571  & 1.213  & 1.022   \\
     & (1.123,1.234) & (1.365,1.555) & (1.472,1.649) & (1.151,1.274) & (0.968,1.076) \\
     \hline
      \rowcolor{LightCyan}
    LLE-MVAR-GH  &  1.182  & 1.452  & 1.586  & 1.214  & 1.023  \\
     \rowcolor{LightCyan}
      & (1.125,1.235) & (1.364,1.563) & (1.468,1.669) & (1.152,1.275) & (0.97,1.076) \\
    \hline
    \hline 
    DM-MVAR-RBF  & 1.387  & 1.534  & 1.902  & 1.560  & 1.131 \\
      & (1.169,6.743) & (1.368,4.634) & (1.609,4.5) & (1.199,9.294) & (0.99,4.655) \\
     \hline 
      \rowcolor{LightCyan}
    LLE-MVAR-RBF  & 1.213  & 1.445  & 1.570  & 1.253  & 1.091  \\
     \rowcolor{LightCyan}
      & (1.135,1.507) & (1.366,1.561) & (1.484,1.657) & (1.159,1.621) & (0.998,1.674) \\
    \hline
  \end{tabular}
  \label{NonLinearManifoldTable}
\end{table}
We also provide the corresponding forecasting statistics as obtained with the proposed ``embed-forecast-lift'' scheme. For embedding in the reduced space, we have taken three (parsimonious) coordinates.
Again, the best performance is obtained with the GHs lifting operator for any combination of manifold learning algorithms and models. Importantly, the  reconstruction errors between the forecasts with the ``full'' MVAR and GPR models trained directly in the original 5D space and the ones obtained with the proposed ``embed-forecast-lift'' scheme are negligible up to a three-digit accuracy for all the five variables. As with the previous case, lifting with RBFs resulted in a poor reconstruction of the high-dimensional space, with forecasting intervals containing the median RMSE value of the naive random walk.

Finally, in Table \ref{modelorder3}, we report the RMSE statistics for the time series produced with the linear model with a model order three (see \eqref{linearmodel3}) as obtained over 100 runs from (a) the naive random walk model applied to the original 5D space, (b) the MVAR models trained in the original 5D data set with model orders one (MVAR(1)) and three (MVAR(3)), and (c) the proposed ``embed-forecast-lift'' method with the embedding applied to the original 5D data set using DM and LLE for embedding with three coordinates. In the reduced-order space we have trained MVAR models with model orders 1 and 3 and used GHs for lifting.
\begin{table}[ht]
 \caption{Linear Stochastic Model with model order three (see \eqref{linearmodel3}). RMSE statistics (median, 5th and 95th percentiles over 100 runs) for 100 simulations, for each of the 5 variables as obtained by: (a) training MVAR(1) and MVAR(3) models in the original 5D feature space, and (b) the proposed scheme with DMs and LLE, MVAR(1) and MVAR(3) models and GHs for lifting. The embedding with DMs and LLE was implemented using three coordinates.}
\centering
 \addtolength{\tabcolsep}{-5pt}
  \begin{tabular}{|l|c|c|c|c|c|}
  \hline
  \rowcolor{Gray}
     Model/Variable  & $y^{(1)}_t$ & $y^{(2)}_t$ & $y^{(3)}_t$ & $y^{(4)}_t$ & $y^{(5)}_t$ \\
    \hline
    \hline
     \rowcolor{LightCyan}
    Random Walk & 2.138  & 2.910 & 1.930  & 3.499  & 2.477  \\
    \rowcolor{LightCyan}
     &(1.898,2.461) &  (2.412,3.551) &  (1.755,2.152) & (2.978,4.195) &  (2.253,2.739) \\
    \hline
    MVAR(1) & 1.629  & 2.121 & 1.382  & 2.567  & 1.840  \\
    & (1.466,1.851) &  (1.817,2.515) &  (1.275,1.517) &  (2.244,2.997) &  (1.699,2.031) \\
    \hline
     \rowcolor{LightCyan}
    MVAR(3) & 1.611  & 2.092  & 1.371  & 2.531 & 1.824 \\
     \rowcolor{LightCyan}
     & (1.449,1.833) &  (1.787,2.483) &  (1.265,1.501) &  (2.206,2.963) & (1.689,2.019) \\
    \hline
    \hline
    DM-MVAR(1)-GH & 1.630  & 2.121  & 1.383  & 2.569 & 1.843  \\ 
     &  (1.467,1.854) &  (1.825,2.515) &  (1.273,1.516) & (2.244,2.998) &  (1.697,2.037) \\ 
     \hline 
   \rowcolor{LightCyan}
    LLE-MVAR(1)-GH & 1.651 & 2.159 & 1.397  & 2.619 & 1.866 \\
    \rowcolor{LightCyan}
     & (1.472,1.905) & (1.832,2.567) & (1.281,1.545) & (2.269,3.080) & (1.705,2.084) \\
     \hline 
     \hline 
    DM-MVAR(3)-GH & 1.621 & 2.103  & 1.376  & 2.546  & 1.835 \\
     & (1.455,1.839) &(1.791,2.484) & (1.267,1.509) & (2.217,2.972) & (1.694,2.028) \\
    \hline
      \rowcolor{LightCyan}
    LLE-MVAR(3)-GH & 1.649 & 2.149 & 1.393 & 2.608 & 1.862  \\
     \rowcolor{LightCyan}
    & (1.473,1.902) & (1.820,2.572) & (1.277,1.550) & (2.252,3.065) & (1.703,2.094) \\
  \hline
  \end{tabular}
  \label{modelorder3}
\end{table}
The best results were obtained when using the proposed scheme with DMs for embedding and a model order three for the training of the MVAR model in the corresponding manifold. Importantly, the implementation of DM-MVAR(3)-GH succeeds in reproducing quite well the results obtained by training MVAR with a maximum delay of three in the original space.

\subsection{FOREX Trading}
Here, we assess the performance of the proposed forecasting framework in the FOREX trading application described earlier, under the annualized Sharpe (SH) ratio (\cite{Sharpe1994}) of the constructed risk parity portfolio. The returns are computed by  Eq. (\ref{eq:RiskParityPortfolio}). The basic formula for calculating the SH ratio is given by:
\begin{equation}
    \text{SH} = \frac{\mu_{\Pi} - R_f}{\sigma_{\Pi}},
\end{equation}
where  $\mu_{\Pi}$ is the sample average returns of the risk parity portfolio and $\sigma_{\Pi}$ the corresponding volatility, while the $R_f$ is the risk-free rate, usually set equal to the annual yield of the US Treasury Bonds. Here, for our analysis, we have set $R_f=0$ that is a fair approximation of the reality. 

The underlying dynamics of the FOREX market is in general non-autonomous, at least over a long time period. In principle, there are time-varying exogenous factors including macroscopic economic indices,  social interaction and mimesis (see e.g. \cite{Papaioannou2013}, where machine learning has been used to forecast FOREX taking into account Twitter messages), but also seasonal factors and rare events (such as the recent COVID19 pandemic), which influence FOREX over time. This comes in contrast to the synthetic time series that we have examined here, and also other autonomous models that have been considered in other studies that served as benchmarks to assess the performance of the various schemes (see also the discussion in the conclusions, Section~\ref{Conclusion}). Hence, to cope with the such changes of the FOREX market and in general of financial assets, one would set up a sliding/rolling window, then train models within the rolling window, and perform (usually) one-day ahead forecasts for trading purposes.

\begin{figure}[ht!]
    \centering
    \hspace*{-1.9cm}
    \includegraphics[scale=0.4]{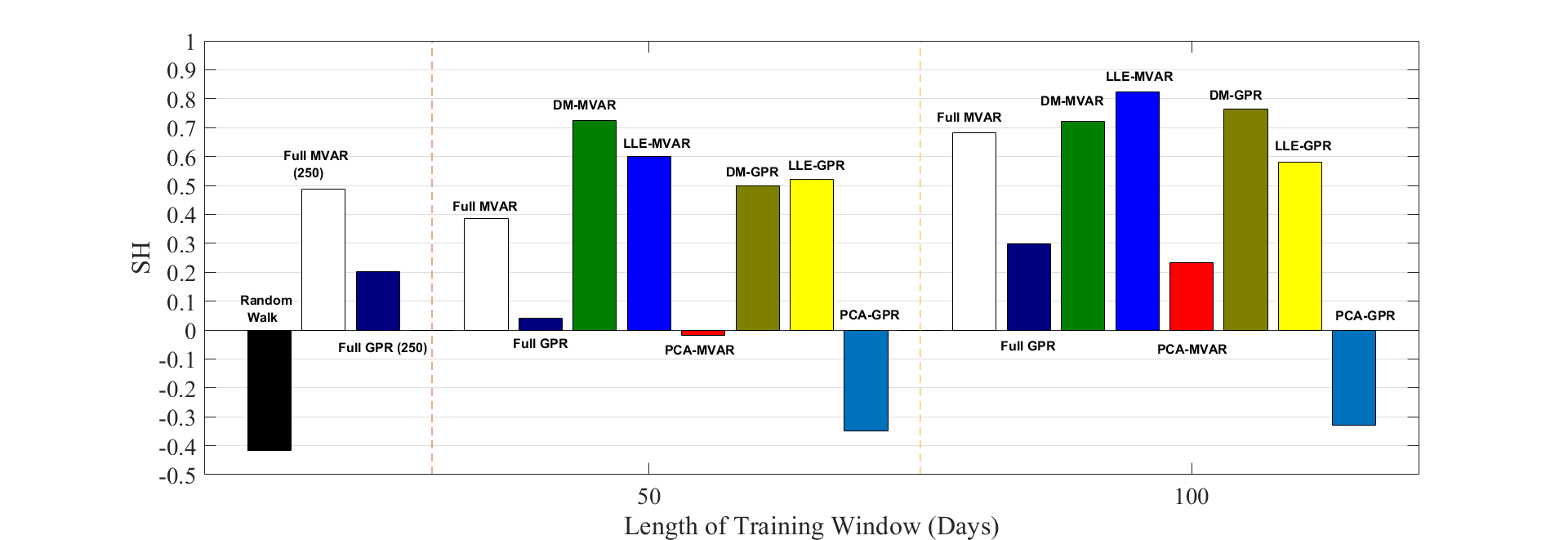}
    \caption{FOREX trading. One-day-ahead predictions. Sharpe Ratios obtained with the proposed framework (using DMs and LLE for embedding, MVAR and GPR for prediction at the embedded space, and GHs for lifting), as well as with PCA, random walk, and MVAR and GPR models implemented directly in the original space. A rolling window of 250 trading days and 3 vectors were used for the embedding, while the MVAR and GPR models in the embedded space were trained using the last 50 or 100 points of the rolling window.}
    \label{fig:SharpeBarPlot_3_AllinOne}
\end{figure}

For our illustrations, we assessed the performance of the proposed scheme based on the so-called risk parity trading strategy, using a 1-year (250 trading days) embedding rolling window and the last 50 or 100 days of the 250-day rolling window for training. The forecasting performance of the proposed scheme using DMs and LLE for embedding with 3 coordinates, MVAR and GPR for prediction, and GHs for lifting was compared against the naive random walk and the full MVAR and GPR models trained and implemented directly in the original space. A comparison against the linear embedding provided by PCA with the same number of principal components was also performed. Figure~\ref{fig:SharpeBarPlot_3_AllinOne} depicts the SH ratios obtained with the various methods. As it is shown, the proposed schemes using the DM and LLE algorithms for embedding outperform all the other schemes, when considering the same size of the training window. In particular, the highest SH ratios ($\sim $0.83) are obtained with the combination of LLE and MVAR within the 100-day training window, followed by the combination of DMs and GPR in the 100-day training window, resulting in an SH ratio $\sim 0.76$. The third ($\sim 0.73$) and fourth  ($\sim 0.72$) larger SH ratios result again from the implementation of the DMs and MVAR within the 50-day and 100-day training windows, respectively. On the other hand, the naive random walk (black bar) resulted in a negative SH ratio ($\sim \; -0.42$). Negative SH ratios (of the order of $-0.35$) resulted also from the implementation of PCA with GPR (yellow bars) for both sizes of the training window, while for the 50-day training window, the combination of PCA with MVAR resulted in an almost zero (but still negative) SH ratio. With PCA, a (small) positive SH value ($\sim 0.23$) was obtained with MVAR within the 100-day training window. The full MVAR and GPR models, trained and implemented directly in the original space, produced positive SH ratios, but still smaller than those of the DMs and LLE schemes. In particular, within the 50-day training window the full MVAR (GPR) model resulted in an SH ratio of $\sim 0.39$ ($\sim 0.04$), while within the 100-day training window the full MVAR (GPR) model resulted in an SH ratio of $\sim 0.68$ ($\sim 0.3$). Finally, we also tested the forecasting performance of the full MVAR and GPR models using the 250-day training window. Within this configuration, the MVAR model yielded an SH ratio of $\sim 0.49$, while the GPR model an SH ratio of $\sim 0.2$.

\section{Conclusions and discussion\label{Conclusion}}

We proposed a numerical framework based on manifold learning for forecasting time series, which is composed of three basic steps: (i) embedding of the original high-dimensional data in a low-dimensional manifold, (ii) construction of reduced-order models and forecasting on the manifold, and (iii) reconstruction of the predictions in the high-dimensional space.

As mentioned in the introduction, the task of forecasting is different from that attempted for the reconstruction of high-dimensional models of dynamical systems in the form of ODEs or PDEs, in four main aspects. First, for real-world data sets, the existence of a relatively smooth low-dimensional manifold is not guaranteed as in the case of well-defined dynamical models (see e.g. the discussion in \cite{gajamannage2019nonlinear}). Second, non-stationary dynamics, which in general are not an issue when dealing with the approximation of dynamical systems, pose a major challenge for a reliable/consistent forecasting. It should be noted that the stationarity assumption may be required to hold true even for interpolation problems. For example, the implementation of the MVAR models, but also Gaussian Processes with a Gaussian kernel, requires the stationary covariance function assumption to be satisfied (see e.g. the discussion in \cite{rasmussen2003gaussian} and \cite{cheng2015time}). Third, when dealing with real-world time series, such as financial time series, the number of available snapshots (even at the intraday frequency of trading) is limited, in contrast to the size of temporal data that one can produce by model simulations. In such cases, the quest for beating the ``curse of dimensionality'' using the correct (parsimonious) embedding and modelling is stronger. Finally, in the case of smooth dynamical systems, as the dynamics emerges from the same physical laws as expressed by the underlying ODEs, PDEs or SDEs, what is usually sought is a single global manifold (a single geometry). However, in many complex problems such as those in finance, the parametrization of the manifold (if it exists) may change over time. Thus, one would seek for a family of (sequential-in-time) submanifolds, which can be ``identified'' within a rolling window approach.

Here, we assessed the performance of the proposed numerical framework by implementing and comparing different combinations of manifold learning algorithms (LLE and DMs), regression models (MVAR and GPR) and reconstruction operators (RBFs and GHs) on various problems, including synthetic stochastic data series and a real-world data set of FOREX time series. By doing so, we first showed that for the synthetic time series, the proposed ``embed-predict-lift'' framework, and in particular the one implementing DMs for embedding and GHs for lifting, reconstructs almost perfectly the predictions obtained in the high-dimensional space with the full regression models. We also showed that using a method such as the LU decomposition for the solution of the RBF interpolation should be in general avoided for lifting, as it might result in an ill-posed problem. However, in this case one could compute the least-norm solution using e.g. the Moore-Penrose pseudoinverse. Moreover, for large scale almost singular linear systems, one can use GMRES and a singular preconditioner as proposed by \cite{elden2012solving}. 
We intend to study the performance of such an approach in a future work.
For the FOREX forecasting problem, we used the standard Sharpe Ratio metric and the parity risk portfolio to assess the forecasting performance of the proposed numerical scheme within a rolling window framework. The proposed scheme with DMs and LLE, MVAR and GPR models and GHs for lifting outperformed all other ``conventional schemes'', i.e. the full MVAR and GPR models implemented directly at the original space as well as the scheme that used PCA for the task of embedding. At this point we should note, that one could use different schemes such as the autoencoders (\cite{kramer1991nonlinear,chen2018molecular}), reservoir computing (\cite{lukovsevivcius2009reservoir,pathak2018model}) and LSTMs (\cite{greff2016lstm,vlachas2018data}). We aim at performing an extensice comparison of the above methodologies in a future work. 
Finally, we note that in order to cope with the generalization property and the topological instability issues (see e.g. \cite{balasubramanian2002isomap}) that arise in implementing kernel-based manifold learning algorithms when the data set is not dense on the manifold, and/or in the presence of ``strong'' stochasticity, one can resort to techniques such as the constant-shift one to appropriately adjust the graph metric, and techniques for the removal of outliers from the data set and the construction of smooth geodesics (\cite{choi2007robust,wang2012geometric,gajamannage2019nonlinear}). We intend to explore the efficiency of such approaches in a future work. 
\section*{Acknowledgements}
This Work partially supported by Gruppo Nazionale per il Calcolo Scientifico - Istituto Nazionale di Alta Matematica (GNCS-INdAM), Italy, and by the Italian program Fondo Integrativo Speciale per la Ricerca (FISR) - B55F20002320001. The work of IGK was partially supported by the US Department of Energy.


\end{document}